# Incentivizing Peer-to-Peer Energy Trading in Microgrids

Amir Noori
*Faculty of Electrical Engineering,*
*K.N. Toosi University of Technology*
Tehran, Iran
anoori@email.kntu.ac.ir

Babak Tavassoli
*Faculty of Electrical Engineering,*
*K.N. Toosi University of Technology*
Tehran, Iran
tavassoli@kntu.ac.ir

Alireza Fereidunian
*Faculty of Electrical Engineering,*
*K.N. Toosi University of Technology*
Tehran, Iran
fereidunian@eetd.kntu.ac.ir

*Abstract*— **Recent trends express the impact of prosumers and small energy resources and storages in distribution systems, due to the increasing uptake of renewable resources. This research studies the effect of coordination of distributed resources with the utility grid and the role of prosumers in the operation of renewable microgrids. We formulated this problem as a social welfare maximization problem followed by employing the dual decomposition method to decompose it into sub-problems of the microgrid, distributed generators, prosumers, and consumers. Moreover, the corresponding power balance mechanism via price adjustment can be viewed as a *Walrasian tatonnement* process. Specifically, prosumers and consumers compete to adjust their energy exchange with other agents to maximize their profit gained by renewable emission reduction benefits while minimizing the associated cost of energy. To this end, we have adopted a peer-to-peer energy trading mechanism based on continuous double auction that can be viewed as a multi-parametric quadratic problem. Finally, we proposed a distributed adaptive algorithm that determines strategies as well as payment and assignment rules. The numerical result suggests that the proposed method can incentivize peer-to-peer energy trading while improving the cost fairness problem and the peak-to-average ratio.**

## I. INTRODUCTION

In power systems, generation must follow the consumption on a real-time basis or even more strict, due to the high electricity storage costs. Renewable energies such as wind and solar characterized by high variability and uncertainty in production has posed several challenges regarding the operation and balancing of power systems [1]. Over the last years, several solutions are proposed to mitigate uncertainty and variability in renewables such as electric, heat, and gas [2], utilizing small and large-scale energy storage systems [3], harnessing the flexibility of demand-side resources [4], or employing better forecasting and coordination algorithms [5-7].

Recently, the integration of renewable energy resources into distribution systems to maintain supply and demand balance has gained more interest [8-9]. One of the main challenges of such methods is how to motivate the independent owners of these resources to behave cooperatively to maximize social welfare. Centralized approaches such as Virtual Power Plants (VPPs) [10-11] are proposed to orchestrate a group of various resources; including distributed generators, storages, and loads to assimilate large generators and loads. However, it does not provide enough motivation for prosumers to contribute to the demand and supply balance. Instead, other market-based coordination methods such as Local Energy and Flexibility Markets (LEFM) [12], Transactive Energy (TE) Systems [13], Peer-to-Peer (P2P) energy trading systems [14-15] try to turn individual consumers from passive to active managers of their networks.

However, these methods need to address several computational difficulties and practical challenges regarding scalability, data security, privacy, regulatory issues, interoperability, to name a few. Currently, to the best of our knowledge, the P2P energy trading paradigm under different available technologies (distributed energy resources (DERs) and consumer level communications and control) is the most viable solution for the distributed energy market. In [16], a centralized P2P electricity trading mechanism based on coalition game theory, formulated as a MILP problem, is proposed in a microgrid with distributed Photo-Voltaics (PVs) and Battery Energy Storage Systems (BESS). In [17], the authors study an energy cost optimization problem among smart homes which are connected for energy sharing and achieved a near-optimal solution based on a bi-level optimization framework. Multiple game-based frameworks, including non-cooperative, evolutionary, and Stackelberg games, are proposed in [18] for P2P energy trading among the prosumers in a community, and two computationally intensive iterative algorithms are proposed to find an equilibrium for the price and seller selection problems. In [19], A peer-to-peer market structure based on a Multi-Bilateral Economic Dispatch (MBED) formulation is introduced and solved based on a Relaxed Consensus + Innovation (RCI) approach in fully decentralized manner.

In this paper, a more realistic problem of energy management in a renewable microgrid with various distributed energy resources are studied. The microgrid is considered in the grid-connected mode that brings up problems for energy pricing and P2P energy trading. Therefore, the impact of P2P energy trading and energy storage systems on some socio-techno-economic issues of microgrids are evaluated. Indeed, it is shown that the amount of excess or shortage of the distributed generated or stored energy in a microgrid changes the market structure from a vertically integrated monopoly to a highly competitive distributed market structure. In general, the problem is more complicated, especially if players behave strategically or other players are waiting to come on. Notably, our analysis does not mean to be comprehensive, and as a preliminary study, we do not yet consider some issues directly like network effects, microgrid islanding, plug-and-play operation, and other stability aspects of microgrids. However, the proposed method seeks to bridge the energy exchange with utility grid and P2P energy trading in a local market, so it enhances the interoperability of the utility grid and microgrids. Moreover, it properly reflects the availability of renewable energy on retail prices. Indeed, it may provide benefits even for consumers without PV panel and energy storage. However, its main advantage is for prosumers who earn more profit from emission reduction benefits that it may also be leveraged as a reward to incentivize P2P energy trading.

The remainder of this paper is organized as follows. In Section II, we introduce the system model. Section III is devoted to formulating the problem. Numerical results are presented in Section IV, and Section V concludes the paper.

## II. SYSTEM MODEL

In this section, we introduce our detailed system model.

### A. Energy System Model

We consider a smart power distribution system that consists of agents with communication and computation capabilities. Specifically, we focus on a typical microgrid equipped with a Microgrid Control system (MCS) and consists of several renewable Distributed Generators (DGs), and a set of prosumers (with PV and/or storage) and consumers, as depicted in Fig. 1. We assume that the MCS, as a rational agent, simultaneously takes on the operator and electricity retailer roles. To this end, the MCS as a retailer predicts the demand and the power generated from local energy resources to determine the net energy demand and then buys the required electricity from the local energy providers and wholesale market. We assume that agents with renewable resources receive an emission reduction benefit. Besides, we assume that end-users (i.e., prosumers and consumers) are equipped with a Home Energy Management System (HEMS) that manages their storage, consumption, and P2P trading patterns. Each prosumer has a preferred energy demand that is partly supplied by its PV and storage. Like consumers, the rest of its required energy is supplied from the utility grid and other prosumers.

We consider a timeframe $\mathcal{T}$ with time granularity $\Delta t$ (i.e., $t_k = k\Delta t \in \mathcal{T}, k \in N$). For the sake of brevity, we also write $t$ instead of $t_k$. We may also discard the indices $t$ in the following. However, the proposed model supports any number of timeslots and time granularities. We also define the communication graph among end-users (prosumers and consumers) by an undirected graph $\mathcal{G}_t = (\mathcal{V}_t, \mathcal{E}_t)$ at the time step $t$ that consists of a set $\mathcal{V}_t$ of vertices and an edge set $\mathcal{E}_t \subseteq \mathcal{V}_t \times \mathcal{V}_t$. Moreover, We assume that $\mathcal{G}_t$ is not fully connected and is a bipartite graph whose vertices compose of two disjoint subsets associated with selling and buying agents.

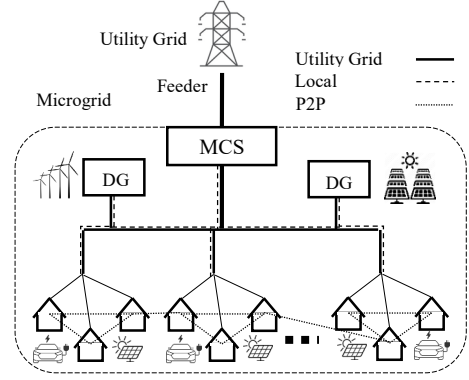

Fig. 1 System Model

At every time step $t$, we assume that each agent $i$ can only communicate with other agents (e.g., $j$) in another subset in its neighborhood, denoted by $j \in \mathcal{N}_i$. However, all agents can communicate with the MCS, while the MCS disseminates retail prices to all.

### B. Distributed Generators

Distributed Generators (DGs) like PV panels and micro-CHPs are the most essential part of a microgrid that maybe used to serve a single load or maybe dedicated to providing clean or reliable power to all customers. Indeed, DGs also can earn revenue by selling back power to the microgrid. The objective of a DG owner $l \in \mathcal{N}_l$ is to maximize the profit given by

$$\max_{p_l} \quad \lambda\, p_l - C_l(p_l) \qquad (1)$$
$$\text{s.t.} \quad \underline{p}_l \leq p_l \leq \overline{p}_l$$

The quadratic cost function is natural candidate among several practical examples

$$C_l(p_l) = a_l\, p_l^2 + b_l\, p_l + c_l$$

Where $p_l$, $\underline{p}_l$ and $\overline{p}_l$ are the power level, the lower and the upper power limits of DG $l \in \mathcal{N}_l$, respectively. $a_l, b_l, c_l$ are appropriate constants known to the DG $l \in \mathcal{N}_l$.

### C. End-users

We consider a set $\mathcal{N}$ of prosumers and consumers. They are the end-user of an energy system. We assume that each end-user $i \in \mathcal{N}$ incurs a cost (or may gain a profit) $\lambda p_i^G$ as a result of purchasing (or selling) $p_i^G$ unit of energy from (or to) the microgrid or other consumers. We assume that each end-user schedules its demand $p_i^d$ a day ahead and also predicts his production $p_i^g$. Therefore, the objective of each end-user is to maximize its profit from selling energy minus the corresponding cost of generating (or purchasing) energy as follows

$$\max_{p_i^G,\, p_{ij},\, p_i^{es}} \lambda p_i^G - C_i(p_i^G) \qquad (2)$$
$$\text{s.t.} \quad p_i^d - p_i^g - p_i^{es} - p_i^G + \sum_{j \in \mathcal{N}_i} p_{ij} = 0$$
$$p_{ij} + p_{ij} = 0 \quad \forall j \in \mathcal{N}_i$$
$$e_{i,t+1}^{es} = e_{i,t}^{es} + \eta p_{i,t}^{es}$$
$$\underline{p}_i^G \leq p_i^G \leq \overline{p}_i^G$$
$$\underline{p}_{ij} \leq p_{ij} \leq \overline{p}_{ij}$$
$$\underline{p}_i^{es} \leq p_i^{es} \leq \overline{p}_i^{es}$$

$$0 \leq e_i^{es} \leq \overline{e}_i^{es}$$

where constraints are the power balance constraint, reciprocity constraints, the energy storage system dynamic, and boundary power constraints, respectively. The power limits are chosen as follow $\underline{p}_{ij}, \overline{p}_{ij} \leq 0$ when delivering power to other consumers, and $\underline{p}_{ij}, \overline{p}_{ij} \geq 0$ otherwise, and $\underline{p}_i^G \cdot \overline{p}_i^G \leq 0 \,(\geq 0)$ when exporting (importing) power to (from) the microgrid, and $\underline{p}_i^{es} \leq 0, \overline{p}_i^{es} \geq 0$. We also assume a quadratic cost function as following

$$C_i(p_i^G) = a_i \,(p_i^G)^2 + b_i \, p_i^G + c_i.$$

### D. Microgrid Control System

Although microgrids can take numerous forms, three typical models are prevalent [20]; DSO monopoly, prosumer consortium, and the free market. In this paper, we adopt a generic prosumer consortium model. In this model, the MGC system is responsible for the overall system operation behind the Point of Common Coupling (PCC). Indeed, this involves several technical and economic issues. However, the technical issues such as dual-mode switching from grid-connected to island mode, power quality, stability and control, and protection are out of the scope of this paper. Here, we only consider the emission reduction benefits cost.

From an economic perspective, the MCS, as a retailer, makes profit by purchasing energy from the utility grid and selling it back to consumers at a higher price to cover its capital and operating costs. However, it can reduce its cost by supplying power from renewable resources. Therefore, the objective of the MCS is given by

$$\max_{p_i^G, p_l, P^G} \lambda P^G + C^{ERB}\left(\sum_i p_i^g \cdot p_l \cdot P^G\right) - C(P^G) \quad (3)$$

$$\text{s.t.} \quad \sum_i p_i^G - \sum_l p_l - P^G = 0$$

where $p_i^G$ and $P^G$ are the power exchange with prosumer $i$ and the utility grid, respectively. The cost function is modeled as quadratic functions given by

$$C(P^G) = a \,(P^G)^2 + b \, P^G + c$$

where $a, b, c$ are appropriate constants known to the MCS, and the Emission Reduction Benefit function are defined as following [21]

$$C^{ERB}\left(\sum_i p_i^g \cdot p_l \cdot P^G\right)$$
$$= 1/T \sum_{i=1}^{4} \left(\left(\sum_i p_i^g + p_l\right) - P^G\right) E_i C_i\right)$$

### III. MARKET CLEARING AND P2P ENERGY TRADING

In this section, we assume that all distributed resources are connected according to the communication graph $\mathcal{G}$ and are also equipped with smart controllers or operated by rational individuals who consider objectives from the previous section. We first formulate the market-clearing problem as a social welfare maximization problem, and then we provide a P2P energy trading mechanism among end-users with a distributed algorithm to coordinate energy pricing and P2P energy trading.

### A. Problem Description

We assume that the MCS gathers an estimate of the net generation and consumption of suppliers and consumers, then calculates the required energy import from the main grid a day ahead and as a retailer purchases the required energy from the wholesale market. In addition we assume that DGs sell their excess generation to the microgrid. However, prosumers can decide on whether sell their excess generation to the microgrid or other consumers. Therefore, the microgrid owner as a retailer clear the market at a retail price $\lambda$ such that profits are maximized and costs are minimized while no demand goes unserved and all local resources are fully utilized. However, it is not always clear how wholesale electric prices translate into retail prices. In perfect competition, the cost-competitive sources of clean energy such as solar, and wind are incorporated as an appropriate source of energy, also known as the merit order effect. In [22], the authors discuss that how the strategic behavior of thermal generators with diverse energy portfolios (including renewable supplies) offsets the price declines due to the merit order effect. They conclude that these strategic players can neutralize the merit order effect when renewable supply is high. However, P2P energy trading can properly reflect the availability of renewable energy on retail prices. The social welfare maximization problem can be described by the following optimization problem

$$\max_{p_i^G, p_j, P^G} -\sum_i C_i(p_i^G) - \sum_l C_l(p_l) - C(P^G) \quad (4)$$
$$+ C^{ERB}\left(\sum_i p_i^G \cdot p_l \cdot P^G\right)$$

$$\text{s.t.} \quad \sum_i p_i^d - \sum_i p_i^g - \sum_i p_i^G - \sum_l p_l = 0$$

with all other constraints of the MCS, prosumers and distributed generators. Note that we removed the term $\sum_i \sum_{j \in \mathcal{N}_i} p_{ij} = 0$.

This problem is a convex problem with linear constraints, and we can solve it efficiently in a centralized manner. However, it is desirable to solve this problem in a decentralized way due to the computational burden and privacy issues. In the next section, we propose Lagrangian dual decomposition method to solve this large-scale problem in a decentralized way.

### B. Structure and Coordination

The corresponding Lagrangian of the problem (4) is defined as following

$$L(\{p_i^G, p_{ij}; \forall i, j \in \mathcal{N}_i\}, \{P_l; \forall l \in \mathcal{N}_l\}, P^G, \lambda)$$
$$= -\sum_i C_i(p_i^G) - \sum_l C_l(p_l)$$
$$- C(P^G)$$
$$+ C^{ERB}\left(\sum_i p_i^G \cdot p_l \cdot P^G\right) \quad (5)$$
$$+ \lambda(\sum_i p_i^d - \sum_i p_i^g - \sum_i p_i^G$$
$$- \sum_l p_l)$$

where $\lambda$ is the dual variable of the power balance constraint. Note that the reduced cost function $C^{ERB}(\sum_i p_i^G. p_l. P^G)$ is decomposable with respect to its arguments. Now, we can decompose the global Lagrangian into several individual Lagrangian terms associated to each agent

$$L_i(p_i^G. \lambda) = -C_i(p_i^G) + C^{ERB}(p_i^G) + \lambda (p_i^d - p_i^g - p_i^G - \sum_{j \in \mathcal{N}_i} p_{ij})$$

$$L_l(p_l. \lambda) = -C_l(p_l) + C^{ERB}(p_l) + \lambda p_l$$

$$L_G(P^G. \lambda) = -C(P^G) + C^{ERB}(P^G)$$

Agent solve their sub-problems in parallel, while the dual variable update provides coordination. However, this algorithm can also be viewed as a form of *walrasian tâtonnement* mechanism or price adjustment process [23]. In a competitive market, the market operator act via price adjustment to reach a market equilibrium by increasing or decreasing the price of a good based on an excess of demand or supply, respectively. Primal update of the MCS and DG owners can be described as follows

$$p_{l.k+1} = \arg min\{ L_l(p_{l.k}. \lambda_k)\} \quad (6)$$

$$p_{k+1}^G = \arg min\{ L(p_k^G. \lambda_k)\} \quad (7)$$

subject to other constraints of DG $l$ and the MCS, while the dual update is given by

$$\lambda_{k+1} = \lambda_k + \rho \left( \sum_i p_i^d - \sum_i p_i^g - \sum_i p_{i.k+1}^G - \sum_j p_{l.k+1} \right) \quad (8)$$

To solve the end-users' sub-problems, it is necessary to develop the P2P energy negotiation mechanism and bidding and offering strategy of prosumers.

*C. P2P Energy Negotiation and Strategies*

Each end-user maximizes its welfare based on proper allocation of power between demand and supply. This allocation includes the power purchased (or supplied) from (to) the grid, the storage system, and P2P trading with neighbors. The end-users optimization problems depend on the exchanged power, which necessitates designing a distributed resource allocation mechanism. Market-based mechanisms, such as double auctions, will represent a widespread means of interaction for P2P energy trading. Three components of these mechanisms; bidding and offering strategies, negotiation schemes, and clearing mechanisms, are introduced in the following.

We assume that users define their asks and bids based on their preferences and the corresponding cost in each time slot $t$. Indeed, they employ an adaptive algorithm to adjust the quantity and the price of the exchanged power. We denote the price and the quantity of P2P power exchanged among end-user $i \in \mathcal{N}$ by $\Theta_i = \{\theta_i = (\lambda_{ij}. p_{ij}) | j \in \mathcal{N}_i\}$. Each end-user $i$ sends its ask/bid to their neighborhoods $j \in \mathcal{N}_i$.

This problem can be viewed as an allocation problem, which is a multi-parametric Quadratic program, given by

$$\max_{p_i^G, \{p_{ij}, \forall j \in \mathcal{N}_i\}} C^{ERB}(p_i^g. p_i^G) - C_i(p_i^G) - \lambda p_i^G - \sum_{j \in \mathcal{N}_i} p_{ij} \lambda_{ij} \quad (9)$$

$$s.t. \quad p_i^d - p_i^g - p_i^{es} - p_i^G - \sum_{j \in \mathcal{N}_i} p_{ij} = 0$$

We solve this problem iteratively in a distributed way. The exchanged power with grid is determined by a gradient descent method with a fixed set of parameters $\Theta_i$ as follows

$$p_{i.k+1}^G = \arg min \ f(p_{i.k}^G)|_{\theta_{ij.k}} \quad (10)$$

where $f$ is the objective function of (9). Then, P2P energy trading update are determined by

$$p_{ij.k+1} = p_{ij.k} + \mu_1(p_{ij.k} + p_{ji.k}) - \mu_2(p_i^d - p_i^g - p_{i.k}^G - \sum_{j \in \mathcal{N}_i} p_{ij.k}) \quad (11)$$

with the following price update rule

$$\lambda_{ij.k+1} = \lambda_{ij.k} - \alpha_{ij}(\lambda_{ij.k} - \lambda_{ji.k}) \quad (12)$$

where $\mu_1$, $\mu_2$ are sequences of positive numbers with persistent excitation, so that the sum of the sequence diverges, and $\alpha_{ij}$ represents the share of each trade from emission reduction benefit given by

$$\alpha_{ij} = \pi_{ij} \frac{\lambda^2(|p_{ij}| + \delta)}{(|p_i^g| + |p_i^G| + \delta)}$$

where $\delta$ is a constant, and $\pi_{ij} = \pi'_{ij} C^{ERB}(p_i^g. p_i^G)$, $\pi'_{ij}$ is the bilateral trading coefficient. The proposed algorithm are listed in Table 1.

*Remark.* To avoid further complexity and keep the problem simple, here we ignored the charging and discharging schedule of the storage system, and we assume a pre-schedule timing based on the availability of renewables, on-peak/off-peak intervals, and without P2P energy trading. Simulations suggest that it would be a good approximation.

Table 1 P2P Energy Trading Negotiations

| Algorithm 1 Double Auction Mechanism |
|---|
| At each time slot $t \in \mathcal{T}$ and for all end-user $i \in \mathcal{N}$ |
| Calculate Initial quantities and prices ($k = 1$) |
| $\lambda_{ij.k} = \lambda \pm \Delta\lambda_{ij}, p_{i.k}^G = p_{i.t}^d - p_{i.t}^g - p_{i.t}^{es}, p_{ij.k} = 1/n \ p_i^G$ |
| While $\Delta p_{i.k}^G > \varepsilon$ do |
|    Send $\hat{p}_{ij.k}, \hat{\lambda}_{ij.k}, (\pi_{ij})$ to neighborhoods $j \in \mathcal{N}_i$ |
|    Compute $p_{i.k+1}^G$ (10) |
|    Update quantity and prices (11) - (12) |
|    IF $(|p_{ij.k+1} + \hat{p}_{ji.k}| \leq \varepsilon_p) \ \& \ (|\lambda_{ij.k+1} - \hat{\lambda}_{ji.k}| \leq \varepsilon_\lambda)$ then |
|      Send Req to $j$ and freeze $\hat{p}_{ij.k}$ and $\hat{\lambda}_{ij.k}$ |
|    IF (received Ack from $j$) then |
|      Establish P2P trade with $j$ with $\hat{p}_{ij.k}$ and $\hat{\lambda}_{ij.k}$ |
|    $k = k + 1$ |
| End |

## IV. NUMERICAL RESULTS

We use numerical simulation to demonstrate the effectiveness of the proposed Algorithm. Here, we consider a full communication graph.

### A. Model implementation and data

The proposed strategy is investigated in simulation on a typical low voltage microgrid test system which is depicted in Fig. 2. This microgrid is in a grid-connected mode and comprises two community-scale PV plant DGs (each 25 kW), 10 households, including 5 with PV and BESS and 5 without PV and BESS. The gross capacity of solar panels and BESSs are assumed to be between 3kW to 5kW, and 3kWh to 7kWh, respectively. The forecasted 24-h load profile, generation of PV panels of prosumers, and the retail energy price update are shown in Fig. 3. The simulations are performed in the MATLAB environment on an Intel core i5 2.50 GHz CPU running Windows 10. We choose a sampling time of one hour. Simulations are performed over one day. The coefficients of cost functions for DGs and prosumers ($a$ and $b$) are chosen between 0.03 to 0.05 c$/kWh$^2$, and 0.2 to 0.5 c$/kWh, respectively. Those are 0.12 c$/kWh$^2$, and 2 c$/kWh for the MCS. The coefficients $E_i$ and $C_i$ for emission reduction benefit function are outlined in Table 2.

Table 2 Environmental Costs of Greenhouse Gases [21]

| Greenhouse gas | $SO_2$ | $NO_X$ | $CO$ | $CO_2$ |
|---|---|---|---|---|
| Emission (kg/kWh) - E | 0.00993 | 0.00646 | 0.00155 | 1.07 |
| Environmental cost ($/kg) - C | 0.97 | 1.29 | 0.16 | 0.0037 |

### B. Simulation and Performance Analysis

In the proposed algorithm, the MCS gather the supply and demand information of all participants and then predicts the net demand of the microgrid and the retail price based on (7), depicted in Fig. 3c. Then, the MCS broadcasts it to DGs and prosumers. Afterwards, DGs and prosumers respond to this price by adjusting their generation level (6) and scheduling their battery energy storage systems and incorporating P2P energy trading (algorithm 1), respectively. However, these depend on forecasts of load profile, PV power output, and market signals including the retail price. Again, MCS updates the retail price based on new forecasts and disseminates the updated prices among end-users. This power adjustment process continues until the retail price converges (depicted in Fig. 3c).

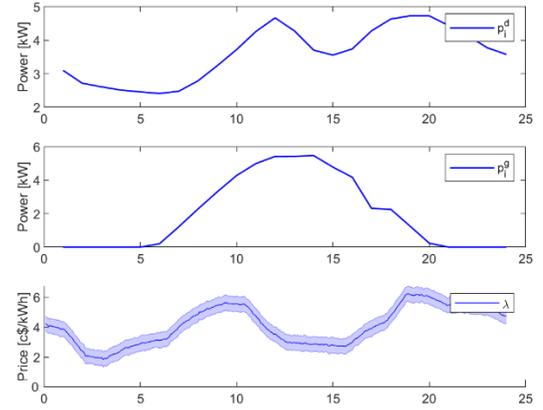

Fig. 3 a) load profile of P1, b) PV power output of P1, and c) the retail price update

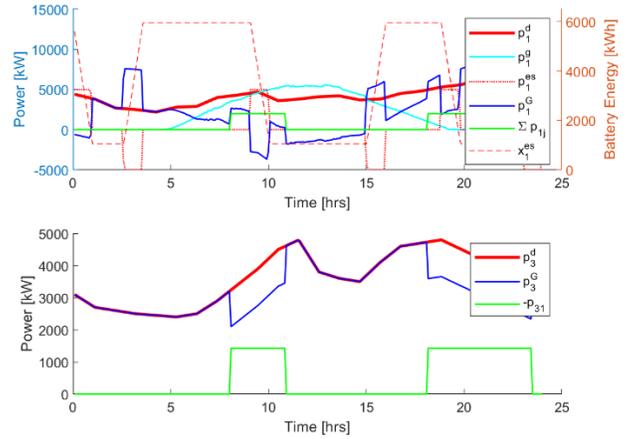

Fig. 4 The power balance of a) prosumer P1 and b) consumer C3

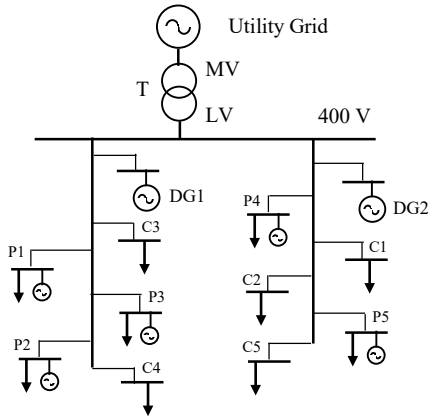

Fig. 2 Test system

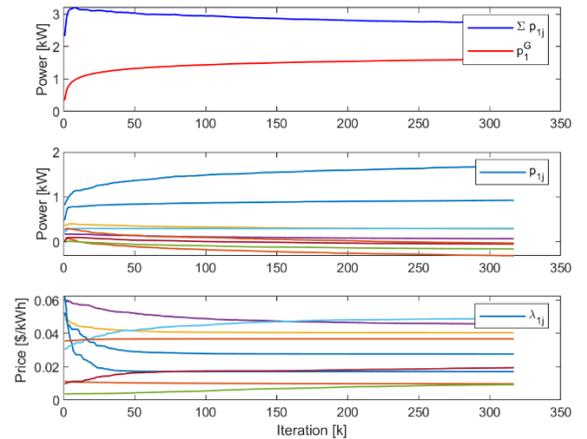

Fig. 5 a) The total P2P energy trading and the exchanged power with utility grid (P1), b) the power of P2P trading, c) the price of P2P trading

In Fig. 4 the Power balance of the prosumer P1 selling energy to the consumer C3 is depicted. The P2P power traded between P1 and C3 is about 500W and P2P price is about 0.04 $/kWh. As depicted in Fig. 4, the proposed P2P energy trading improve the cost fairness problem among the participating peers and it reduces the Peak-to-Average Ratio (PAR). In Fig. 5, the result of the proposed P2P energy trading mechanism in time-slot $t = 10$ is depicted. End-users aim to maximize their welfare by adjusting the ratio of power supplied from the grid and that from P2P trading. In addition, prosumers empowered by the emission reduction benefit can adjust the P2P price through the bilateral trading coefficient $\pi_{ij}$.

## V. Conclusion

In this paper, we proposed a coordination mechanism to deal with the challenge brought by the high penetration of renewable energy resources in distribution systems. We bring into account the local generation, energy storage systems, and P2P energy trading. We employed a distributed double auction mechanism within the price adjustment process to simultaneously clear the local market and enable P2P energy trading. The numerical result suggests that the proposed method improves the interoperability of the microgrid and distributed energy resources. Future studies will include the investigation of network effects and uncertainties.